\title{ ~~\\ On primes in arithmetic progression
having a prescribed primitive root.~II}
\author{Pieter Moree}
\documentclass[12pt]{article}
\usepackage{amssymb,latexsym,amsfonts}
\def\@ptsize{2}
\newtheorem{Thm}{Theorem}
\newtheorem{lem}{Lemma}

\newtheorem{Prop}{Proposition}
\newcommand{\qed}{\hfill $\Box$}

\begin{document}
\date{}
\maketitle
\centerline{\small \it Dedicated to Prof.W. Narkiewicz,on the occasion of his 70th birthday}
\begin{abstract} 
\noindent
Let $a$ and $f$ be coprime positive integers. Let $g$ be 
an integer.
Under the Generalized Riemann Hypothesis (GRH) it follows by 
a result of H.W. Lenstra that the set of primes $p$ such
that $p\equiv a({\rm mod~}f)$ and $g$ is a primitive root
modulo $p$ has a natural density.
In this note this density is explicitly evaluated with an Euler product
as result.
\end{abstract}
\section{Introduction}
Let $g\in \mathbb Z\backslash \{0\}$. 
Let ${\cal P}_g$ denote the
set of primes $p$ such that $g$ is a primitve root mod $p$. 
It was conjectured by Emil
Artin in 1927 that ${\cal P}_g$ is infinite in case $g$ is in ${\cal G}$, the set
of integers not equal to $-1$ or a square. 
Clearly, if $g$ is not in ${\cal G}$ then ${\cal P}_g$ is finite.
There is no integer $g$ for
which {\it Artin's primitive root conjecture} (as it is usually called) has been proved. However, Heath-Brown 
\cite{HB} in a classical paper, basing himself on a breakthrough paper of
Gupta and Murty \cite{GM},
established a result which implies, for example, that there are at most
two primes $q$ for which ${\cal P}_q$ is finite. In 1967 Hooley
\cite{Hooley}
established Artin's conjecture under the assumption of
the Generalized Riemann Hypothesis (GRH). Moreover, he showed
that under that assumption the set ${\cal P}_g$ has a natural density,
which he evaluated (his result is
Theorem \ref{main} below with $a=1$ and $f=1$).
It turns out that this density is a rational
number (depending on $g$) times the Artin constant $A$,
with 
$$A=\prod_p\left(1-{1\over p(p-1)}\right)=0.3739558136 1920228805 4728054346 \cdots$$
(Here and in the rest of the paper the notation $p$ is used to indicate primes.)\\
\indent In connection with his study of Euclidian number fields, Lenstra
\cite{L}
considered the distribution over 
arithmetic progressions of the primes in ${\cal P}_g$.
Let ${\cal P}_{a,f,g}$ denote the set of primes such that $g$ is a primitive
root mod $p$ and $p\equiv a({\rm mod~}f)$. {}From Lenstra's work it follows
that, under GRH,  ${\cal P}_{a,f,g}$ has a natural density. 
\begin{Thm} {\rm \cite{L}}.
\label{lens}
Put $\zeta_m=e^{2\pi i/m}$. Let $f\ge 1$ and $1\le a\le f,~(a,f)=1$. Let
$\sigma_a$ be the automorphism of $\mathbb Q(\zeta_f)$ determined by
$\sigma_a(\zeta_f)=\zeta_f^a$. Let $c_a(n)$ be $1$ if the restriction
of $\sigma_a$ to the field $\mathbb Q(\zeta_f)\cap \mathbb Q(\zeta_n,g^{1/n})$
is the identity and $c_a(n)=0$ otherwise. Put
$$\delta(a,f,g)=\sum_{n=1}^{\infty}{\mu(n)c_a(n)\over 
[\mathbb Q(\zeta_f,\zeta_n,g^{1/n}):\mathbb Q]}.$$
Then, assuming GRH, we have
$$\pi_{g}(x;f,a)=\delta(a,f,g){x\over \log x}+
O_{f,g}\left({x\log \log x\over \log^2x}\right),$$
where $\pi_{g}(x;f,a)$ denotes the number of primes $p\le x$ that
are in ${\cal P}_{a,f,g}$.
\end{Thm}
In the light of the apparent arithmetical complexity of Lenstra's formula, the
following 
relatively simple expression for $\delta(a,f,g)$ may come as a bit of a surprise.
\begin{Thm} 
\label{main}
Let $g\in {\cal G}$.
Let $h\ge 1$ be the largest integer such that $g$ is an $h$th power.
Let $\Delta$ denote the discriminant of the quadratic 
field $\mathbb Q(\sqrt{g})$. Let $f\ge 1$ and $1\le a\le f,~(a,f)=1$.
Let $b=\Delta/(f,\Delta)$. Put
$$\gamma=\cases{(-1)^{b-1\over 2}(f,\Delta) & if $b$ is odd;\cr
1 & otherwise.}$$
Put
$$A(a,f,h)=\prod_{p|(a-1,f)}(1-{1\over p})
\prod_{p\nmid f\atop p|h}(1-{1\over p-1})\prod_{p\nmid f\atop p\nmid h}
\left(1-{1\over p(p-1)}\right)$$
if $(a-1,f,h)=1$ and $A(a,f,h)=0$ otherwise. Then
$$\delta(a,f,g)={A(a,f,h)\over \varphi(f)}\left(1+
({\gamma\over a}){\mu(2|b|)\over 
\prod_{p|b,~p|h}(p-2)\prod_{p|b,~p\nmid h}(p^2-p-1)}\right).$$
Here $({\cdot\over \cdot})$ denotes the Kronecker symbol.
\end{Thm}
On writing $g=g_1g_2^2$, with $g_1$ squarefree and the $g_i$ integers,
we see that $\Delta=g_1$
if $g_1\equiv 1({\rm mod~}4)$ and $\Delta=4g_1$ otherwise. Note that
$b$ is odd if and only if $g_1\equiv 1({\rm mod~}4)$
or $g_1\equiv 2({\rm mod~}4)$ and $8|f$
or $g_1\equiv 3({\rm mod~}4)$ and $4|f$. Theorem
\ref{main} can be easily reformulated in terms of $g_1$.
\begin{Thm} 
\label{main1}
Let $a,f,g,h,\gamma$ and $A(a,f,h)$ be
as in Theorem {\rm \ref{main}}. Let 
$$\beta={g_1\over (g_1,f)}{\rm ~and~}\gamma_1=\cases{(-1)^{\beta-1\over 2}(f,g_1) & if $\beta$ is odd;\cr
1 & otherwise.}$$ We have
$$\delta(a,f.g)={A(a,f,h)\over \varphi(f)}\left(1-
({\gamma_1\over a}){\mu(|\beta|)\over 
\prod_{p|\beta,~p|h}(p-2)\prod_{p|\beta,~p\nmid h}(p^2-p-1)}\right)$$
in case $g_1\equiv 1({\rm mod~}4)$ or $g_1\equiv 2({\rm mod~}4)$ and
$8|f$ or $g_1\equiv 3({\rm mod~}4)$ and $4|f$ and
$$\delta(a,f,g)={A(a,f,h)\over \varphi(f)},$$
otherwise. 
\end{Thm}
{}From Theorem \ref{main} various known results can be rather easily deduced.
That is the subject of 
Section \ref{applic}. The remaining sections
are devoted to proving Theorem \ref{main}.\\

\noindent {\tt Remark}. The results in this paper were first proved in 
1998 in preprint MPIM1998-57 \cite{primoud}. Not much later Lenstra et al. \cite{LMS}, for a preview
see Stevenhagen \cite{S}, found a method which allows one to give a more
conceptual and elegant proof of Theorem 2 and hence I
abstained from trying to publish the present paper. However, at present
it does not seem clear whether \cite{LMS} will be ever completed. Since
also the method presented here requires a rather more modest input
from algebraic number theory than \cite{LMS} and is quite different, I
now decided to publish (an updated and polished version of) MPIM1998-57.\\
\indent  The proof of the main result
of part I of this series, Theorem 4 of \cite{pie3}, depends crucially on Theorem \ref{main}
(see part 7 of \S 5 for a brief description of Theorem 4 of \cite{pie3} and its proof).

\section{Some facts from algebraic number theory}
\label{algfacts}
Although Theorem \ref{lens} suggests differently, the problem 
of evaluating $\delta(a,f,g)$ really only
involves cyclotomic fields, quadratic subfields
and their composita. In this section some facts 
concerning these fields relevant for the proof of
Theorem \ref{main} are discussed. We start by recalling some properties
of the Kronecker symbol, a rarely covered topic in books on number theory (but
see e.g. Cohen \cite[pp. 36-39]{Cohen}).\\
\indent The Kronecker symbol $({a\over b})$ is the extension of the Jacobi
symbol to $(\mathbb Z\backslash \{0\})^2$ obtained by setting $({a\over -1})={\rm sign}(a)$
and $({a\over 2})=({2\over a})$ for $a$ odd ($({a\over 2})=0$ for $a$ even), and extending
by multiplicativity. All symbols of the form $({\cdot\over \cdot})$ in this paper will be Kronecker symbols.
For two nonzero integers $m$ and $n$ write
$m=2^{\nu_2(m)}m_1$ and $n=2^{\nu_2(n)}n_1$, with $m_1$ and $n_1$ odd. Then
$$({n\over m})=(-1)^{((m_1-1)(n_1-1)+({\rm sign}(m)-1){\rm sign}(n)-1))/4}
({m\over n}).$$
This is the law of quadratic reciprocity formulated in terms of the Kronecker symbol. Furthermore,
if $D$ is a discriminant then (see Cohen \cite[Theorem 2.2.9]{Cohen})
\begin{equation}
\label{reci2}
({D\over a+kD})=({D\over a}).
\end{equation}
\indent The following
lemma allows one to determine all quadratic subfields of a given cyclotomic
field (for a proof see e.g. \cite[p. 163]{W}).
\begin{lem}
\label{classical}
Let $\mathbb Q(\sqrt{d})$ be a quadratic
field of discriminant $\Delta_d$.
Then
the smallest cyclotomic field containing $\mathbb Q(\sqrt{d})$ is
$\mathbb Q(\zeta_{|\Delta_d|})$.
\end{lem}
Consider the cyclotomic field $\mathbb Q(\zeta_f)$. There are $\varphi(f)$
distinct automorphisms determined uniquely by $\sigma_a(\zeta_f)=\zeta_f^a$, with
$1\le a\le f$ and $(a,f)=1$. We need
to know when the restriction of such an automorphism to a given
quadratic subfield of $\mathbb Q(\zeta_f)$ is the identity. 
\begin{lem} 
\label{twee}
Let $\mathbb Q(\sqrt{d})\subseteq \mathbb Q(\zeta_f)$ be a quadratic
field of discriminant $\Delta_d$.
We have $\sigma_a|_{\mathbb Q(\sqrt{d})}=$id if and only if
$({\Delta_d\over a})=1$.
\end{lem}
{\it Proof}. Notice that, by Lemma \ref{classical}, we can restrict to
the case where $f=|\Delta_d|$. 
Define $\chi$ by $\chi(a)=\sigma_a(\sqrt{d})/\sqrt{d}$, 
$1\le a\le |\Delta_d|,~(a,\Delta_d)=1$.
Then $\chi$ is the uniqe non-trivial character of the character group 
of $\mathbb Q(\sqrt{d})$. As is well-known (see e.g. 
\cite[p. 437]{N1}) the primitive character
induced by this is $({\Delta_d\over a})$. Using Lemma \ref{classical}, we
see that $\chi$ is a primitive character mod $|\Delta_d|$. Thus
$\chi(a)=({\Delta_d\over a})$. Now $\sigma_a|_{\mathbb Q(\sqrt{d})}=$id if and
only if $\chi(a)=({\Delta_d\over a})=1$. \qed \\

\noindent {\tt Remark}. It is also possible to prove 
Lemma \ref{twee} using quadratic reciprocity and properties
of Gauss sums.\\
\indent The next result can be proved by a trivial generalization of an
argument given by Hooley \cite[pp. 213-214]{Hooley}.
\begin{lem}
\label{cycdegree}
Let $g\in {\cal G}$ and let $h$ be the largest positive integer such that $g$ is
an $h$th power. Let $\Delta$ be the discriminant of the quadratic field $\mathbb Q(\sqrt{g})$.
Suppose that $k|r$ and $k$ is squarefree. Put $k_1=k/(k,h)$ and
$n(k,r)=[\mathbb Q(\zeta_r,g^{1/k}):\mathbb Q].$ Then
\begin{itemize}
\item[i)] if $k$ is odd, $n(k,r)=k_1\varphi(r);$
\item[ii)] if $k$ is even and $\Delta\nmid r,$ $n(k,r)=k_1\varphi(r);$
\item[iii)] if $k$ is even and $\Delta|r,$ $n(k,r)=k_1\varphi(r)/2.$
\end{itemize}
\end{lem}
{\tt Remark}. Without the condition that $k$ be squarefree the latter lemma becomes
much more complicated to state, 
see e.g. Moree \cite[Lemma 1]{modern}. Indeed, in the general case different definitions of $h$ and $\Delta$ 
appear to be the natural ones: one writes $g=\pm g_0^h$ with $g_0>0$ and $h$ as large as possible and
denotes with $\Delta$ the discriminant of $\mathbb Q(\sqrt{g_0})$. Partly because of these
different definitions it is not so straightforward to check that the general lemma implies  Lemma
\ref{cycdegree}. This is left as an exercise for the interested reader.\\
\indent The next lemma together with Lemma \ref{twee} allows one to compute $c_a(n)$.
\begin{lem} 
\label{gemier}
Let $g\in {\cal G}$. Let $\Delta$ denote
the discriminant of $\mathbb Q(\sqrt{g})$. Let $n\ge 1$ be squarefree and
$f\ge 1$ be arbitrary. Put $b=\Delta/(f,\Delta)$. Put
$$\gamma=\cases{(-1)^{b-1\over 2}(f,\Delta) & if $b$ is odd;\cr
1 & otherwise.}$$
Then
$$\mathbb Q(\zeta_f) \cap \mathbb Q(\zeta_n,g^{1/n})=\cases{\mathbb Q(\zeta_{(f,n)},
\sqrt{\gamma}) &
if $n$ is even, $\Delta\nmid n$ and $\Delta|{\rm lcm}(f,n)$;\cr
\mathbb Q(\zeta_{(f,n)}) & otherwise.}$$
\end{lem}
{\it Proof}. On noting that  $\varphi((f,n))\varphi({\rm lcm}(f,n))=
\varphi(f)\varphi(n)$, it easily follows, using Lemma \ref{cycdegree},
that 
\begin{equation}
\label{graad}
[\mathbb Q(\zeta_f)\cap \mathbb Q(\zeta_n,g^{1/n}):\mathbb Q(\zeta_{(f,n)})]=2
\end{equation}
if $n$ is even, $\Delta\nmid n$ and $\Delta|{\rm lcm}(f,n)$ and
$\mathbb Q(\zeta_f) \cap \mathbb Q(\zeta_n,g^{1/n})=\mathbb Q(\zeta_{(f,n)})$
otherwise. Thus we may assume that  $n$ is even, $\Delta\nmid n$ and $\Delta|{\rm lcm}(f,n)$. Notice that 
this implies that $b$ is odd.
It is easy to check that $\gamma$ is a fundamental
discriminant. Since clearly $\gamma|f$ it follows by Lemma \ref{classical} that
$\sqrt{\gamma}\in \mathbb Q(\zeta_f)$. Note that $\Delta|{\rm lcm}(f,n)$ implies that
${\Delta\over (f,\Delta)}|{n\over (f,n)}|n$. From this it follows that $(f,\Delta)|(f,n)$
would imply $\Delta|n$, contradicting the assumption $\Delta\nmid n$. This contradiction shows
that $\sqrt{\gamma}\not\in \mathbb Q(\zeta_{(f,n)})$. Note that $\Delta/\gamma$ is a fundamental
discriminant. Since $b|n$ it follows by Lemma \ref{classical} that
$\sqrt{\Delta/\gamma}\in \mathbb Q(\zeta_n)$ and thus, since $n$ is even, 
$\sqrt{\gamma}\in \mathbb Q(\zeta_n,g^{1/n})$.
In summary we have that
$\sqrt{\gamma}\not\in \mathbb Q(\zeta_{(f,n)})$, 
$\sqrt{\gamma}\in \mathbb Q(\zeta_f)$ and 
$\sqrt{\gamma}\in \mathbb Q(\zeta_n,g^{1/n})$. These inclusion relations
together with (\ref{graad}) yield the result. \qed

\section{Euler products}
In this section we prove some
results that will help us to write down the Euler
product of the sums encountered
in the proof of Theorem \ref{main}.
\begin{Prop}
\label{multiplicative}
Let $f,h\ge 1$ be integers. Then the function $w:\mathbb N\rightarrow \mathbb N$
defined by 
$$w(k)={k\varphi({\rm lcm}(k,f))\over (k,h)\varphi(f)}$$
is multiplicative. Furthermore,
\begin{itemize}
\item[i)] if $p\nmid h$ and $p\nmid f$, then $w(p)=p(p-1)$
\item[ii)] if $p\nmid h$ and $p|f$, then $w(p)=p$
\item[iii)] if $p|h$ and $p\nmid f$, then  $w(p)=p-1$
\item[iv)] if $p|h$ and $p|f$, then $w(p)=1$
\item[v)] if $h$ is odd, then $w(2)=2$.
\end{itemize}
\end{Prop}
{\it Proof}. For every multiplicative function $g$ and arbitrary integers $a,b\ge 1$, we obviously
have $g(a)g(b)=g((a,b))g({\rm lcm}(a,b))$. Thus, it is enough to show that $\varphi((k,f))$ is
a multiplicative function in $k$, which is obvious. The remaining part of the result follows
on direct calculation. \qed\\

The multiplicativity of $w$ plays an important role in the proof of the following lemma.
\begin{lem} 
\label{uiltje}
Let $a,f,h\ge 1$ be integers
with $1\le a\le f,~(a,f)=1$ and $h$ odd. Let
$\Delta$ be a discriminant of a quadratic number field. Let
$b=\Delta/(f,\Delta)$. Put
$$S(b)=\sum_{n=1,~\Delta|{\rm lcm}(n,f)\atop a\equiv 1({\rm mod~}(f,n))}
^{\infty}{\mu(n)\over w(n)}.$$
Let $S_2(b)$ denote the same sum as $S(b)$ but with the restriction that
$2|n$. 
Then 
$$S(b)=-{\mu(2|b|)A(a,f,h)\over \prod_{p|b}(w(p)-1)}.$$
Furthermore, $S_2(b)=-S(b)$.
\end{lem}
{\it Proof}. If $b$ is even, then the summation in $S(b)$ runs over
non squarefree $n$ only and hence $S(b)=0$. Next assume that $b$ is odd.
We have
\begin{eqnarray}
S(b) &=& \sum_{n=1,~b|n/(f,n)\atop a\equiv 1({\rm mod~}(f,n))}^{\infty}
{\mu(n)\over w(n)}=\sum_{d|(a-1,f)}\sum_{n=1,~(f,n)=d\atop b|n/d}
^{\infty}{\mu(n)\over
w(n)}\nonumber\\
&=& \sum_{d|(a-1,f)}{\mu(d)\over w(d)}\sum_{n=1,~(f,n)=1\atop b|n}^{\infty}
{\mu(n)\over w(n)}={\mu(|b|)\over w(|b|)}\sum_{d|(a-1,f)}{\mu(d)\over w(d)}
\sum_{n=1,~(f,nb)=1\atop (b,n)=1}^{\infty}{\mu(n)\over w(n)}.\nonumber
\end{eqnarray}
Now by assumption $\Delta$ is a discriminant and $b$ is odd. This implies
that $(f,b)=1$ and $b$ is squarefree. Thus
\begin{eqnarray}
S(b) &=& {\mu(|b|)\over w(|b|)}\prod_{p|(a-1,f)}(1-{1\over w(p)})
\prod_{p\nmid fb}(1-{1\over w(p)})\nonumber\\
&=& {\mu(|b|)\over w(|b|)}
\prod_{p|(a-1,f)}(1-{1\over w(p)})
\prod_{p\nmid f}(1-{1\over w(p)})\prod_{p|b}(1-{1\over w(p)})^{-1}\nonumber\\
&=& {\mu(|b|)A(a,f,h)\over \prod_{p|b}(w(p)-1)}=-{\mu(2|b|)A(a,f,h)\over \prod_{p|b}(w(p)-1)},\nonumber 
\end{eqnarray}
where we used that $w(p)>1$ for $p|b$ and
$$\prod_{p|(a-1,f)}(1-{1\over w(p)})\prod_{p\nmid f}(1-{1\over w(p)})=A(a,f,h),$$ an identity immediately obtained on invoking
Proposition \ref{multiplicative}.\\
\indent On using that $w(2)=2$ along with the fact that both $\mu$ and $w$ are multiplicative, 
we find that
$$S_2(b)=-{1\over 2}\sum_{2\nmid n,~\Delta|{\rm lcm}(2n,f)\atop a\equiv 1({\rm mod~}(f,2n))}
{\mu(n)\over w(n)}.$$
First assume that $f$ is even. Then $a$ is odd and hence, for $n$ odd, 
$a\equiv 1({\rm mod~}(f,2n))$ if and only if $a\equiv 1({\rm mod~}(f,n))$.
It follows that
$$S_2(b)=-{1\over 2}\sum_{2\nmid n,~\Delta|{\rm lcm}(n,f)\atop a\equiv 1({\rm mod~}(f,n))}
{\mu(n)\over w(n)}.$$ One easily checks that the latter formula also holds if $f$ is odd.
On noting that
$$S_2(b)+\sum_{2\nmid n,~\Delta|{\rm lcm}(n,f)\atop a\equiv 1({\rm mod~}(f,n))}
{\mu(n)\over w(n)}=S(b),$$
we infer that $S_2(b)=-S(b)$. \qed

\section{Proof of Theorem \ref{main}}
\noindent i) The case 
$b$ is odd and $({\gamma\over a})=1$. Using Lemma \ref{gemier}, Lemma
\ref{twee}
and the observation above, we find that
$c_a(n)=1$ in case $a\equiv 1({\rm mod~}(f,n))$
and $c_a(n)=0$ otherwise. This together with Theorem \ref{lens} and
Lemma \ref{cycdegree} implies 
that 
\begin{eqnarray}
\label{star}
\varphi(f)\delta(a,f,g)&=&\sum_{n=1\atop 2\nmid n}^{\infty}{\mu(n)\over w(n)}
+\sum_{n=1,~2|n\atop \Delta\nmid {\rm lcm}(n,f)}^{\infty}{\mu(n)\over w(n)}+
2\sum_{n=1,~2|n\atop \Delta|{\rm lcm}(n,f)}^{\infty}{\mu(n)\over w(n)}.\nonumber\\
&=&\sum_{n=1}^{\infty}{\mu(n)\over w(n)}+
\sum_{n=1,~2|n\atop \Delta|{\rm lcm}(n,f)}^{\infty}{\mu(n)\over w(n)}
={\rm I}_1+S_2(b),
\end{eqnarray}
where furthermore in each sum we restrict to those 
integers $n$ such that
$a\equiv 1({\rm mod~}(f,n))$.\\ 

\noindent ii) The case $b$ is odd and $({\gamma\over a})=-1$. Using Lemma \ref{gemier},
Lemma \ref{twee}, the observation in the beginning of this proof 
and (\ref{star}), we find
$$\varphi(f)\delta(a,f,g)={\rm I}_1+S_2(b)-2\sum_{{2|n,~\Delta\nmid n\atop
\Delta|{\rm lcm}(f,n)}\atop a\equiv 1({\rm mod~}(f,n))}
{\mu(n)\over w(n)}.$$
Now
$$\sum_{2|n,~\Delta|{\rm lcm}(f,n)\atop
a\equiv 1({\rm mod~}(f,n))}{\mu(n)\over w(n)}=
\sum_{{2|n,~\Delta\nmid n\atop \Delta|{\rm lcm}(f,n)}\atop
a\equiv 1({\rm mod~}(f,n))}{\mu(n)\over w(n)}
+\sum_{2|n,~\Delta|n\atop a\equiv 1({\rm mod~}(f,n))}{\mu(n)\over w(n)}.$$
In case $\Delta\equiv 0({\rm mod~}4)$, the latter sum is obviously zero.
In case $\Delta\equiv 1({\rm mod~}4)$, a necessary condition for the latter
sum to be non-zero is that $a\equiv 1({\rm mod~}(f,2\Delta))$. 
By (\ref{reci2}) it then follows that $({\gamma\over a})=({\gamma\over 1})=1$, this
contradiction shows that the latter sum is always zero. It thus follows that
$\varphi(f)\delta(a,f,g)={\rm I}_1-S_2(b).$ \\

\noindent iii) The case $b$ is even. In this case we cannot have that $n$ is squarefree, 
$2|n$, $\Delta\nmid n$ and $\Delta|{\rm lcm}(f,n)$ (for this would imply that $b$ is odd)
and as in case i) we find that $c_a(n)=1$ in case $a\equiv 1({\rm mod~}(f,n))$
and $c_a(n)=0$ otherwise. As in case i) we then find that
$\varphi(f)\delta(a,f,g)=I_1+S_2(b)$.\\

Note that in all three cases we have $\varphi(f)\delta(a,f,g)=I_1+({\gamma\over a})S_2(b)$.
Now, on using Proposition \ref{multiplicative}, 
\begin{eqnarray}
{\rm I}_1&=&\sum_{d|(a-1,f)}\sum_{(f,n)=d}{\mu(n)\over w(n)}=\sum_{d|(a-1,f)}{\mu(d)\over w(d)}\sum_{(f,n)=1}{\mu(n)\over w(n)}\nonumber\\
&=& \prod_{p|(a-1,f)}(1-{1\over w(p)})\prod_{p\nmid f}(1-{1\over w(p)})
=A(a,f,h).\nonumber
\end{eqnarray}
On invoking Lemma \ref{uiltje} it follows that
$$\varphi(f)\delta(a,f,g)=I_1+({\gamma\over a})S_2(b)=A(a,f,h)\Big(1+({\gamma\over a}){\mu(2|b|)\over \prod_{p|b}(w(p)-1)}\Big).$$
On working out the product
using  Proposition \ref{multiplicative} the proof is then completed. \qed

\section{Applications}
\label{applic}
{\bf 1)}. {\it Hooley's theorem}. Setting $a=1$ and $f=1$ in Theorem \ref{main1} we 
obtain Hooley's theorem \cite{Hooley}. \\
{\bf 2)}. $a=1$. Setting $a=1$ in Theorem \ref{main} we
obtain Theorem 4 of \cite{pie2}. Notice
that of all the progressions mod $f$, the
progression $1({\rm mod~}f)$ is the easiest to deal with, since
trivially $c_1(n)=1$ for every $n$.\\
{\bf 3)}. {\it Rodier's conjecture}. Rodier \cite{rod}, in
connection with a coding theoretical problem, conjectured that
the density of the primes in ${\cal P}_2$ such that $p\equiv -1,3$
or $19 ({\rm mod~}28)$ is $A/4$. {}From Theorem \ref{main} it follows
however that, under GRH, this density is $21A/82$ (cf. \cite{rodier}): for
each progression $a_i({\rm mod~}28)$ under consideration we find that
$\delta(a_i,28,2)=A(a_i,28,1)/\varphi(28)=21A/246$..\\
{\bf 4)}. {\it Zero density}. Lenstra \cite[Theorem 8.3]{L}
gave a sketch of a proof of the following result.
\begin{Thm}
Let $g\in {\cal G}$. Then $\delta(a,f,g)=0$ if and only if one of the following holds
\begin{itemize}
\item[i)] $(a-1,f,h)>1$;
\item[ii)] $\Delta|f$ and $({\Delta\over a})=1$;
\item[iii)] $\Delta|3f$, $3|\Delta$, $3|h$ and $({-\Delta/3\over a})=-1$.
\end{itemize}
\end{Thm}
This result very easily follows from Theorem \ref{main}. We leave it to
the reader to show that if $\delta(a,f,g)=0$, then actually ${\cal P}_{a,f,g}$
is finite. In each of these cases, there are obstructions not going beyond quadratic
reciprocity. So, assuming GRH, loosely speaking $\delta(a,f,g)=0$ if and only if there is an elementary
obstruction, or obstruction coming from quadratic reciprocity.\\
{\bf 5)}. {\it Equidistribution}. If $S$ is any set of integers, denote by $S(x)$ the number of integers in $S$
not exceeding $x.$ For given integers $a$ and $f,$ denote by $S(x;f,a)$
the number of integers in $S$ not exceeding $x$ that are congruent to $a$
modulo $f.$ We say that $S$ is weakly uniformly distributed mod $f$ 
(or WUD mod $f$ for short) if
$S(x)\rightarrow \infty$ $(x\rightarrow \infty)$ and for every $a$ coprime
to $f,$
$$S(x;f,a)\sim {S(x)\over \varphi(f)}.$$
The progressions $a({\rm mod~}f)$ such that the latter asymptotic equivalence
holds and $S(x)\rightarrow \infty$ $(x\rightarrow \infty)$ are said to get their {\it fair share} of primes from 
$S$. Thus 
$S$ is {\it weakly uniformly distributed} mod $f$ if and only if all the progressions
mod $f$ get their fair share of primes from $S.$
Narkiewicz \cite{N} has written a nice survey on the state of knowledge
regarding the
(weak) uniform distribution of many important arithmetical sequences.
Let ${\cal D}_g$ denote the
set of natural numbers $f$ such that ${\cal P}_g$ is weakly uniformly
distributed modulo $f$.
\begin{Thm}
\label{equidis}
{\rm (GRH)}.
Let $g\in {\cal G}$ and let
$h$ the largest integer such that $g$ is an $h$th
power. Write $g=g_1g_2^2$ with $g_1$ squarefree.
Assume that not both $g_1=21$ and $(h,21)=7,$ then 
${\cal D}_g,$
the set of natural numbers $f$ such that the set of primes $p$
such that $g$ is a primitive root mod $p$ is weakly uniformly distributed
mod $d,$ equals
\begin{itemize}
\item[i)] $\{2^n:n\ge 0\}$ if $g_1\equiv 1({\rm mod~4});$ 
\item[ii)] $\{1,2,4\}$ if $g_1\equiv 2({\rm mod~4});$ 
\item[iii)] $\{1,2\}$ if $g_1\equiv 3({\rm mod~4}).$ 
\end{itemize}
In the
remaining case $g_1=21$ and $(h,21)=7,$ 
we have ${\cal D}_g=\{2^m3^n:n,~m\ge 0\}.$
\end{Thm}
Using only a formula for $\delta(1,f,g)$ 
and Theorem \ref{lens} in some special cases, this result was first deduced
in \cite{pie2}. Here a shorter proof, using the full force of Theorem
\ref{main1}, is given.\\

\noindent {\it Proof of Theorem} \ref{equidis}. Put
$S_f:=\{A(a,f,h)~|~1\le a\le f,~(a,f)=1\}$. Let $\gamma$ be as in Theorem \ref{main}.
Let us first consider the case where $f=2^m$ for some $m\ge 0$. 
Notice that $|S_{2^m}|=1$. 
By Theorem \ref{main1} we now infer that 
${\cal P}_g$ is WUD mod $2^m$
if
$g_1\equiv 2({\rm mod~}4)$ and $m\le 2$, or $g_1\equiv 3({\rm mod~}4)$
and $m\le 1$. In the remaining case we have
$$({\gamma_1\over a})=\cases{({1\over a}) & if $g_1\equiv 1({\rm mod~}4)$;\cr
({-1\over a}) & if $g_1\equiv 3({\rm mod~}4)$;\cr
({(-1)^{g_1-2\over 4}2\over a}) & if $g_1\equiv 2({\rm mod~}4)$,}$$
and ${\cal P}_g$ is WUD $2^m$ if and only if $({\gamma_1\over a})$ is
trivial character, that is if and only if $g_1\equiv 1({\rm mod~}4)$.\\
\indent It remains to deal with the case where $f$ has an odd prime divisor. First let us
consider the case where $f=q$ with $q$ an odd prime. In case $\beta$ is
even
$$\varphi(f)\delta(1,f,g)=A(1,f,h)\ne A(2,f,h)=\delta(2,f,g)\varphi(f),$$ and we do not have equidistribution. Next
assume that $\beta$ is odd.
If $q=3$, then a short
calculation shows that a necessary and sufficient
condition for equidistribution
to occur, is that $3|g_1,~(3,h)=1,~\mu(|\beta|)=-1$, $g_1\equiv 1({\rm mod~}4)$
and the equation 
$\prod_{p|\beta,~p|h}(p-2)\prod_{p|\beta,~p\nmid h}(p^2-p-1)=5$ has a solution
with $\beta$ is odd. Now notice that $g$ is solution to this if and only
if $g_1=21$ and $(h,21)=7$. Call such a $g$ exceptional. 
If $q\ge 5$, then there exists
$2\le a_1\le q-1$ such that $({\gamma_1\over a_1})=1$ and so, by Theorem \ref{main}, 
$\delta(1,f,g)\neq \delta(a_1,f,g)$.
{}From this it immediately follows that in case
$f$ has an odd prime divisor and
$g$ not exceptional, equidistribution does not occur, for if ${\cal P}_g$ is
WUD mod $f$, then ${\cal P}_g$ must also be WUD mod $\delta$ for every
divisor $\delta$ of $f$. It remains to show that when $g$ is
exceptional we have that ${\cal P}_g$ is WUD
mod $2^m3^n$, with $m\ge 0$ and $n\ge 1$ arbitrary. This follows
easily from Theorem \ref{main} and the calculation done in the
case $f=3$. \qed\\
\noindent {\bf 6)}. {\it Optimal normal basis}. Let $q$ be an odd given prime power and $m$ a natural number. One can wonder whether
there exists an extension ${\mathbb F}_{q^n}$ of ${\mathbb F}_{q^m}$ such that
${\mathbb F}_{q^n}$ has an optimal normal basis over $\mathbb F_q$. A number
theoretic question that arises in this context, see
\cite{CGPT},  is whether there is a prime $p$
such that $p\equiv 1({\rm mod~}m)$ and $q$ is a primitive root modulo $p$ and if
yes to provide a small upper bound for the smallest such prime. The latter part
of the question seems very difficult, but the results obtained in this paper allow
one to shed some light on the first part of the question \cite{GMP}.\\
{\bf 7)}. {\it Asymptotically exact heuristics}. On invoking Hooley's theorem \cite{Hooley} it can be shown, 
see \cite{primsom, pie3}, that, under GRH, we have
$$\pi_g(x;2,1)=2\sum_{2<p\le x,~(g/p)=-1\atop (p-1,h)=1}{\varphi(p-1)\over p-1}+
O_{g}\Big({x\log \log x\over \log^2x}\Big),$$
were the sum is not so difficult to evaluate (with $(g/p)$ the
Legendre symbol, $g\in {\cal P}_g$ and $h$ the largest integer such that $g$ is an $h$th
power). One would then expect that
\begin{equation}
\label{bloder}
\pi_g(x;f,a)=2\sum_{2<p\le x,~(g/p)=-1\atop p\equiv a({\rm mod~}f),~(p-1,h)=1}{\varphi(p-1)\over p-1}+
O_{f,g}\Big({x\log \log x\over \log^2x}\Big).
\end{equation}
Unconditionally the sum appearing here can be evaluated with error $O(x\log^{-c}x)$ with $c>0$
arbitrary, see \cite[Theorem 1]{pie3}. It turns out that the main term equals $\delta(a,f,g)$Li$(x)$,
with $\delta(a,f,g)$
as given in Theorem 3. From this the truth of (\ref{bloder}), on GRH, then follows. This
is the main result (Theorem 4) of \cite{pie3}.\\

\noindent {\bf Acknowledgement}. The author thanks F. Lemmermeyer and P. Stevenhagen 
for some helpful remarks.

\medskip\noindent {\footnotesize Max-Planck-Institut f\"ur Mathematik,
Vivatsgasse 7, 53111 Bonn, Germany.\\ 
(e-mail: {\tt moree@mpim-bonn.mpg.de})}

\end{document}